\documentclass{gtpart}     
\title{On phantom maps into co-H-spaces}

\author{James Schwass}
\givenname{James}
\surname{Schwass}
\address{Department of Mathematics\\ Grand Valley State University\newline{} Allendale, MI 49001}
\email{james.schwass@gvsu.edu}
\urladdr{}


\arxivreference{1512.00353}
\arxivpassword{xwep8}

\volumenumber{}
\issuenumber{}
\publicationyear{}
\papernumber{}
\startpage{}
\endpage{}
\doi{}
\MR{}
\Zbl{}
\received{}
\revised{}
\accepted{}
\published{}
\publishedonline{}
\proposed{}
\seconded{}
\corresponding{}
\editor{}
\version{}

\newtheorem*{BSThm}{Bott-Samelson Theorem}      

\usepackage{amsmath,amsfonts,amsthm,amssymb,verbatim,enumerate}
\usepackage[all]{xy}
\usepackage{setspace}
\usepackage{Tabbing}
\usepackage{fancyhdr}
\usepackage{lastpage}
\usepackage{extramarks}
\usepackage{chngpage}
\usepackage{soul,color}
\usepackage{graphicx,float,wrapfig}
\usepackage[toc,page]{appendix}

\newcommand{\hsand}{\hspace{3pc}\mbox{and}\hspace{3pc}}

\renewcommand{\tilde}{\widetilde}
\renewcommand{\bar}{\overline}

\newcommand{\amin}{A^{\min}}
\newcommand{\bmax}{B^{\max}}
\newcommand{\bamin}{\bar{A}^{\min}}
\newcommand{\bbmax}{\bar{B}^{\max}}

\newcommand{\CoH}{\mathbf{CoH}}

\newtheoremstyle{HWProblem}
{1pc}{1pc}{\upshape}{}{\bfseries}{.}{.5em}{}

\theoremstyle{BookTheorem}
\newtheorem{thm}{Theorem} [section]
\newtheorem*{thm*}{Theorem}
\newtheorem{lem}[thm]{Lemma}
\newtheorem{prop}[thm]{Proposition}
\newtheorem{cor}[thm]{Corollary}

\theoremstyle{HWProblem}

\theoremstyle{definition}

\newtheorem{ex}[thm]{Example}

\theoremstyle{rem}

\newtheorem{que}[thm]{Question}

\newcommand{\id}{\mathrm{id}}

\renewcommand{\phi}{\varphi}

\renewcommand{\P}{\mathcal{P}}
 \newcommand{\ZZ}    {\mathbb{Z}}

    \newcommand{\NN}    {\mathbb{N}}
    \newcommand{\FF}    {\mathbb{F}}
    \newcommand{\QQ}    {\mathbb{Q}}

    \newcommand{\om}    {\Omega}
    \newcommand{\s}     {\Sigma}

\newcommand{\Top}{\mathbf{Top}}

\newcommand{\too}{\longrightarrow}
\newcommand{\lra}[1]{\stackrel{#1 }{\longrightarrow}}

\DeclareMathOperator{\limone}{lim^1}

        \DeclareMathOperator{\cat}   {cat}

    \DeclareMathOperator{\Ext}   {Ext}

        \DeclareMathOperator{\conn}{conn}

 \DeclareMathOperator{\Ph}    {Ph}

\newcommand{\LS}{Lusternik-Schnirelmann }

\begin{document}

\begin{abstract}We study the existence of essential phantom maps into co-H-spaces, motivated by Iriye's observation that every suspension space $Y$ of finite type with $H_i(Y;\QQ)\neq 0$ for some $i>1$ is the target of essential phantom maps. We show that Iriye's observation can be extended to the collection of nilpotent, finite type co-H-spaces.
This work hinges on an enhanced understanding of the connections between homotopy decompositions of looped co-H-spaces and coalgebra decompositions of tensor algebras due to Grbi\`{c}, Theriault, and Wu. 
\end{abstract}

\maketitle


\section{Introduction}

We will work in the category $\Top$ of spaces having the homotopy type of a pointed CW complex and pointed maps between them. We will restrict our attention throughout to simply-connected spaces, or their loop spaces. 
A map $X\to Y$ is called a \textbf{phantom map} if for every $n$ the composite
\[
X_n\to X \to Y
\]
is nullhomotopic, where $X_n\to X$ is an $n$-skeleton for some CW structure of $X$. We offer an alternative characterization of this concept to illustrate that the choice of a CW for structure $X$ is insignificant; according to \cite{BousfieldKan}, $X\to Y$ is phantom if and only if $X\to Y \to Y^{(n)}$ is nullhomotopic for every $n$, where $Y^{(n)}$ denotes the $n$th Postnikov approximation of $Y$. 

From the definition and characterization given above, it is clear that a phantom map must induce the zero map on homotopy groups, and on any homology theory, and so these maps appear trivial upon passage to such common algebraic models for topological spaces. On the other hand, phantom maps can be of genuine topological interest. 
The theory of phantom maps has been used by 
Harper and Roitberg \cite{HarperRoitberg}, and Gray \cite{graySNT}, among many others, 
to produce and study examples of distinct homotopy classes of spaces $X$ and $Y$ which have the same $n$-type, i.e. $X^{(n)}\simeq Y^{(n)}$, for all $n$. 
Roitberg \cite{RoitbergAut} has also used the theory of phantom maps to compute the homotopy automorphism groups of particular spaces; in general the computation of homotopy automorphism groups is intractable. These examples serve to illustrate that phantom maps play a significant role in $\Top$. But, since these maps vanish under many of our favorite functors, they prove difficult to study, or even to locate. The purpose of this work is to locate new examples of phantom maps; the analysis of particular invariants of these phantom maps and the structure of the collection of phantom maps will take place elsewhere. 

  The constant map is an obvious example of a phantom map. Of more interest are essential (i.e. homotopically nontrivial) phantom maps, which abound in $\Top$. We offer, as evidence of this fact, the following theorems of Iriye, and McGibbon and M\o ller. We will say a space $X$ is of finite type (over $\ZZ$) if each  $H_n(X;\ZZ)$ and $\pi_n(X)$ is a finitely generated group. 
  We write $\Ph(X,Y)$ for the subset of $[X,Y]$ consisting of homotopy classes of phantom maps. 


\begin{thm}{\rm \cite{iriyePhantomsToSuspensions}\label{thmIriye}}\qua
Suppose $Y\simeq \s X$ is a nilpotent suspension space of finite type. If $H_i(Y;\QQ)\neq 0$ for some $i>1$ then $Y$ is the target of essential phantom maps from finite type domains. 
\end{thm}

\begin{thm}{\rm\cite{McGMoller}}\qua
If $X$ and $Y$ are of finite type and $\Ph(X,Y)$ is not the one point set, then $\Ph(X,Y)$ is uncountably large. 
\end{thm}


In many senses, the concept of a co-H-space is a mild generalization of that of a suspension space. As such, many statements that hold true for the collection of suspension spaces are also true for the collection of co-H-spaces. 
We wondered if one could replace the suspension space $Y$ in Theorem \ref{thmIriye} with any nilpotent co-H-space of finite type. 
Our main result is a positive answer to this question.

%

\begin{thm}\label{thmMain}
Suppose $Y$ is a nilpotent co-H-space with $H_i(Y;\QQ)\neq 0$ for some $i>1$. Then $Y$ is the target of essential phantom maps from finite type domains. 
\end{thm}

The proof of Theorem \ref{thmMain} is comprised of several pieces. For a co-H-space whose rational homology is ``large'' we develop decomposition methods in phantom map theory and appeal to recently-developed highly-structured decompositions of the loopspace of a co-H-space due to Selick, Grbi\`{c}, Theriault, and Wu. For a co-H-space with ``small'' rational homology we exploit strong connections between phantom map theory and rational homotopy theory discovered by McGibbon and Roitberg.

Through the theory of \LS category, this work can be viewed as providing a solution to the case $n=1$ of the following question. Our exposition of \LS category here will  be limited to the following three observations:  $\cat(X)$ is a non-negative integer assigned to a space $X$ which we think of as a measure of the complexity of $X$;  $\cat(X)=0$ if and only if $X$ is contractible;   the spaces of \LS category one are precisely the noncontractible co-H-spaces. 

\begin{que}\label{quePhantomsToFiniteCat} 
Suppose $Y$ has finite type, and $\cat(Y)=n<\infty$. If $H_i(Y;\QQ)\neq 0$ for some $i>1$, is $Y$ the target of essential phantom maps from finite type domains? 
\end{que}

In Section \ref{secPhantomPrelim} we lay out the preliminaries on phantom map theory. 
In Section \ref{secDecompCoH} we describe recently developed connections between coalgebra decompositions of tensor algebras and homotopy decompositions of looped co-H-spaces. In Section \ref{secDecompMethods} we develop techniques to bridge the gap between the decompositions of Section \ref{secDecompCoH} and the theory of phantom maps. Section \ref{secCoH} contains the proof of Theorem \ref{thmMain}. 
Examples and applications are given in Section \ref{secExamples}. 


\paragraph{Acknowledgements}
We would like to thank Jeff Strom, the author's dissertation advisor, for many helpful conversations regarding the content and preparation of this paper. This work grew out of the author's doctoral dissertation at Western Michigan University. We are indebted to Kouyemon Iriye for supplying a critical insight into the proof of Proposition \ref{prop43}, which helped this work reach its maturity.  
We would also like to thank the referee for suggesting the inclusion of additional examples.

\section{Preliminaries}

Localization will play a central role in what is to follow. We assume familiarity with the rudiments of localization; a detailed reference is \cite{MayPonto}.  Since a rationally nontrivial $p$-local space is not of finite type over $\ZZ$, we will have a need for a $p$-local analog of the notion of a finite type space; a space $X$ is of \textbf{finite type over $\ZZ_{(p)}$} if each $H_n(X;\ZZ)$ and $\pi_n(X)$ is a finitely generated $\ZZ_{(p)}$-module. We should note that a space of finite type over $\ZZ_{(p)}$ is necessarily $p$-local. 
Though we will be primarily interested in phantom maps between finite type spaces, we will have occasion to examine phantom maps from finite type domains into targets having finite type over $\ZZ_{(p)}$. 

 \subsection{Background on Phantom Maps} \label{secPhantomPrelim}

%
%

In Section \ref{secTowerPerspective} we describe a critical identification of $\Ph(X,Y)$ with a particular functor which factors through the category of towers of groups. In Section \ref{secPhantomRatEq} we describe connections between phantom map theory and rational homotopy theory that are indispensable in discovering new examples of phantom maps from old, among other things. Most of the material in this section can be found in the wonderful survey article \cite{handbook} of McGibbon.

%

\subsubsection{The Tower Perspective}\label{secTowerPerspective}

By a tower $\{G_n\}$ of groups we mean a diagram 
\begin{equation}\label{diagTower}
\dots \lra {p_{n+1}} G_n \lra {p_n} \dots \lra{p_3} G_2 \lra {p_2} G_1 
\end{equation} 
in the category of groups. We mean something similar by a tower of Abelian groups, or a tower of sets, or really a tower of any sort of gadget - these are $\NN^\mathrm{op}$-shaped diagrams in various categories. A morphism of towers is a natural transformation of $\NN^\mathrm{op}$ shaped diagrams. By $\lim G_n$ we mean the limit of the diagram \eqref{diagTower} in the appropriate category. 

We now set about describing the functor $\limone$. On the category of towers of Abelian groups, by $\limone$ we mean the first derived functor of $\lim$; more concretely, if $\{G_n\}$ is a tower of Abelian groups, then 
 $\lim G_n$ is the kernel and $\limone G_n$ is the cokernel of the map 
\[
\prod G_n \lra{\id-(p_n)} \prod G_n
\]
given by 
\[
(a_1,a_2,\dots)\mapsto (a_1-p_2(a_2),a_2-p_3(a_3),\dots).
\]

Bousfield and Kan \cite[pgs 254--255]{BousfieldKan} extend the definition of $\limone$ to the category of towers of arbitrary groups as follows: 
Given a tower $\{G_n\}$ of groups let $\prod G_n$ act on $\prod G_n$ by 
\[
(g_n)\cdot (x_n) = (g_n x_n (p_{n+1}(g_{n+1})^{-1})),
\]
where $G_{n+1}\lra{p_{n+1}} G_n$ is the structure map in the tower $\{G_n\}$. 
Then $\limone G_n$ is the orbit space of this action. This is important to us because we will have occasion to refer to $\limone G_n$ where $\{G_n\}$ is a tower of not necessarily Abelian groups.

%
%

In particular, if $X$ and $Y$ have the homotopy type of CW complexes, then 
a CW structure for $X$ gives rise to a tower 
$\{[\s X_n,Y]\}$ of (generally non-Abelian) groups; dually the Postnikov tower for $Y$ gives rise to a tower 
$\{[X,\om Y^{(n)}]\}$ of (generally non-Abelian) groups. 
 We now arrive at a fundamental identification in phantom map theory. 

%
%

\begin{cor}{\rm \cite{BousfieldKan}} \qua \label{corPhLimone}
For spaces $X$ and $Y$ there are bijections of pointed sets
\[
\limone[ \s X_n, Y] \cong \Ph(X,Y)\cong \limone [ X , \om Y^{(n)}]. 
\]
\end{cor}  

The identification made in Corollary \ref{corPhLimone} allows for the introduction of algebraic methods for characterizing the condition $\Ph(X,Y)=\ast$.
Given a tower of gadgets (groups, sets, etc.) $\{G_n\}$ let $G_k^{(n)}$ denote the image in $G_k$ of the composite of the structure maps 
\[
G_n\to G_{n-1} \to \dots \to G_k
\]
when $n\geq k$ and for $n< k$ set $G_k^{(n)}=1$. 
This defines, for each $k\geq 1$ a subtower $\{G_k^{(n)}\}$, indexed by $n$, of the tower $\{G_n\}$. Notice that for fixed $k$ the sequence of images $G_k^{(n)}$ are nested; we say the tower $\{G_n\}$ \textbf{satisfies the Mittag-Leffler condition} if all of the nested sequences $G_k^{(n)}$ satisfy a descending chain condition: explicitly, for each $k$ there is some $N$ so that for all $n\geq N$ one has $G_k^{(n)}=G_k^{(N)}$. 

It is well-known that if a tower $\{G_n\}$ satisfies the Mittag-Leffler condition, then $\limone G_n=\ast$. When the tower $\{G_n\}$ is comprised of countable groups, the converse of this statement is also true:

\begin{thm}{\rm \label{limoneCtble}\cite{McGMoller}}\qua 
Suppose $G_n$ is a tower of countable groups. Then $\limone G_n=\ast$ if and only if the tower $G_n$ satisfies the Mittag-Leffler condition. Moreover, if \newline{} $\limone G_n\neq \ast$, then $\limone G_n$ is uncountable large. 
\end{thm}

It is worthwhile to note that when $X$ and $Y$ are of finite type over $\ZZ$ or $\ZZ_{(p)}$ for some prime $p$, then for each $n$ the groups
\[
[\s X_n,Y] \hsand [X,\om Y^{(n)}]
\]
are countable. Theorem \ref{limoneCtble}  will be used to develop decomposition methods in phantom map theory in Section \ref{secDecompMethods}.

\subsubsection{Phantom Maps and Rational Equivalences} \label{secPhantomRatEq}


 McGibbon and Roitberg have characterized the finite type spaces that are not the targets of essential phantom maps from finite type domains in terms of the existence of particular rational equivalences.

\begin{thm}{\rm\cite{McGRoit}\label{thmNoPhantoms}}\qua For a nilpotent, finite type space $Y$, the following are equivalent 
\begin{enumerate}[(i)]
\item  $\Ph(X,Y)=\ast$ for all finite type domains $X$, 

\item $\Ph(K(\ZZ,m),Y)=\ast$ for all $m$, and 
 
   \item  there is a rational equivalence $\prod_\alpha K(\ZZ,m_\alpha)\to\om Y$.  
   \end{enumerate}
\end{thm}

We should note that the direction of the rational equivalence in Theorem \ref{thmNoPhantoms} part (iii) is significant; for any space $Y$ there is a rational equivalence $\om Y \to \prod K(\ZZ,m_\beta)$. 

We will need a $p$-local version of the implication (i) $\Rightarrow$ (iii) of Theorem \ref{thmNoPhantoms}, which we record as Proposition \ref{propPLocalNoPhantoms}. This will be used to establish a lemma in Section \ref{secDecompMethods} required to develop decomposition methods in phantom map theory. 

 We have previously observed that if $X$ and $Y$ are of finite type over $\ZZ$ or $\ZZ_{(p)}$, then the groups 
 \[
 [ \s X_n, Y] \hsand [X, \om Y^{(n)}]
  \] 
  are countable for all $n$. As such, Theorem \ref{limoneCtble} can be used to characterize the condition $\Ph(X,Y)=\ast$ in terms of the Mittag-Leffler condition. This is the main point required to complete the construction of the rational equivalence $\prod K(\ZZ,m_\beta)\to \om Y$ as given by McGibbon and Roitberg, given the hypothesis $\Ph(X,Y)=\ast$ for all finite type domains $X$, and so we arrive at the following partial refinement of Theorem \ref{thmNoPhantoms}.


\begin{prop}\label{propPLocalNoPhantoms}
Suppose $Y$ is nilpotent and has finite type over $\ZZ_{(p)}$. If $\Ph(X,Y)=\ast$ for all finite type domains $X$, then there is a rational equivalence 
\[
\prod K(\ZZ,m_\beta)\to \om Y.
\]
\end{prop}

The converse of this statement could feasibly hold, but we have not yet had occasion to check this. Indeed, if conjugacy classes in $[X,\om Y^{(n)}]$ are of finite cardinality for every $n$, then the converse of Proposition \ref{propPLocalNoPhantoms} can be established using the proof of Theorem \ref{thmNoPhantoms} given by McGibbon and Roitberg in \cite{McGRoit}.

Theorem \ref{thmNoPhantoms} only begins to hint at the connections between phantom map theory and rational homotopy theory. The next result is another glimpse of these strong connections.
We should note that the result stated here is slightly stronger than in \cite{McGRoit}, though the authors' argument establishes the result in light of the observation that $[X,\om Y^{(n)}]$ is a countable group when $X$ and $Y$ are of finite type over $\ZZ$ or $\ZZ_{(p)}$. Before stating the result, we remark that $\Ph(X,Y)$ is a contravariant functor in $X$ and a covariant functor in $Y$.

\begin{thm}{\rm\cite{McGRoit}\label{thmRatSurj}}\qua
Suppose $Y$ and $Y'$ are of finite type over $\ZZ$ or $\ZZ_{(p)}$. If $Y\to Y'$ induces a surjection on $\pi_\ast\otimes\QQ$, then for every finite type domain $X$ the induced map 
\[
\Ph(X,Y)\to \Ph(X,Y')
\]
is surjective. 
\end{thm}

Note that for each prime $p$ and each nilpotent space $Y$ the $p$-localization $Y\to Y_{(p)}$ is a rational equivalence, hence induces surjections on $\pi_\ast\otimes\QQ$, and so we arrive at a corollary which has been well-known in the phantom map literature, and will be one of our primary tools for detecting essential phantom maps.

\begin{cor}\label{corPhantomsLiftThroughLoc}
Suppose $Y$ is a nilpotent, finite type space. If $Y_{(p)}$ is the target of essential phantom maps from finite type domains, then so is $Y$. 
\end{cor}


\subsection{Homotopy Decompositions of Looped Co-H-Spaces} \label{secDecompCoH}

Our jumping off point is the generalized Bott-Samelson theorem, due to Berstein. 

\begin{BSThm}{\rm\cite{BersteinBottSamelson} }\qua
If $Y$ is a simply-connected co-H-space, then there is a natural algebra isomorphism 
\[
H_\ast(\om Y) \cong T(\s^{-1} \tilde{H}_\ast(Y)),
\]
where $H_\ast(\om Y)$ is equipped with the Pontryagin product. Here homology has coefficients in a PID $k$ and $ \tilde{H}_\ast(Y)$ is a free $k$-module.
\end{BSThm}

For the rest of this section we fix a prime $p$; the ground ring for all algebraic objects will be $\FF_p$, the field with $p$ elements. 
 All homology in this section has $\FF_p$ coefficients.  Many of the results of this section remain true if we replace $\FF_p$ with an arbitrary field, though we will have no need for such generality. 
We write $T$ for the free graded tensor algebra functor taking the category of vector spaces to the category of graded algebras.

In the 1980s, F. Cohen, Moore, and Neisendorfer developed a technique fueled by the Bott-Samelson theorem which they use to determine the homotopy exponents of odd dimensional spheres; the difficulty of drawing concrete conclusions regarding homotopy groups of spheres is well-documented, and illustrates the power of this technique. We now loosely outline one component of this program. 
 Cohen, Moore, and Neisendorfer sought out algebraic decompositions of $T(\s^{-1} \tilde{H}_\ast(Y))$, and showed that these algebraic decompositions have geometric realizations in the form of homotopy decompositions of $\om Y$ for $Y=S^{2n+1}$, among a few other specific spaces.

 In \cite{SelickWu} Selick and Wu begin developing functorial analogs of the ad hoc decomposition methods of Cohen, Moore, and Neisendorfer, apparently motivated by the power of these methods, along with a conjecture of F. Cohen. The functorial decomposition methods reach maturity in \cite{DecompOmCoH}, after contributions by Grbi\`{c}, Theriault, Selick, and Wu spanning the course of about a decade.  Before describing these functorial analogs, we lay out some nomenclature and conventions.

Of course as vector spaces $T(V)=\bigoplus_{n\geq 0} V^{\otimes n}$, where $V^{\otimes 0}=\FF_p$. This identifies $V$ as a submodule of $T(V)$. The algebra $T(V)$ is equipped with a unit $\FF_p\to T(V)$ and augmentation $T(V)\to \FF_p$ defined by inclusion of and projection onto  $\FF_p=V^{\otimes 0}$, respectively. 
The tensor algebra $T(V)$ is naturally endowed with the structure of a Hopf algebra by declaring the elements of $V$ to be primitive. More explicitly, since $T(V)$ is the free algebra on $V$, the linear map $V\to T(V)\otimes T(V)$ given by $v\mapsto 1\otimes v + v\otimes 1$ extends uniquely to a map of algebras $\Delta:T(V)\to T(V)\otimes T(V)$, giving a comultiplication on $T(V)$. One can check that the unit and augmentation are morphisms of coalgebras and algebras, respectively, and so we have given $T(V)$ the structure of a Hopf algebra. 
This discussion serves to illustrate that we can think of the tensor algebra functor $T$ as taking its values in the categories of algebras, coalgebras, or Hopf algebras. We will specify which category we mean to take for the target of the functor $T$ if there is potential for confusion.

 A \textbf{natural coalgebra retract} of $T$ is a functor $A$ from vector spaces to coalgebras equipped with natural transformations $A\lra I T$ and $T\lra R A$ so that $RI$ is the identity natural transformation on $A$. 
 A \textbf{natural coalgebra decomposition} of $T$ is a pair of functors $A,B$ from vector spaces to coalgebras equipped with natural coalgebra isomorphisms $T\cong A\otimes B$. Since $\otimes$ is the categorical product in the category of coalgebras, which happens to be a pointed category, it follows that if $T\cong A\otimes B$ is a natural coaglebra decomposition, then both $A$ and $B$ are natural coalgebra retracts of $T$. 
 A \textbf{natural sub-Hopf algebra} of $T$ is a subfunctor $B$ from vector spaces to Hopf algebras. A natural sub-Hopf algebra $B$ of $T$ is \textbf{coalgebra split} if $B$ is a natural coalgebra retract of $T$ when regarded as a functor into the category of coalgebras. 

 We will write $\CoH_{(p)}$ for the category of $p$-local co-H-spaces and co-H-maps between them. 
 A \textbf{natural homotopy retract} of $\om:\CoH_{(p)}\to \Top$ is a functor $\bar A:\CoH_{(p)}\to \Top$ equipped with natural transformations $\bar A\lra I \om$ and $\om\lra R \bar A$ so that $RI$ is naturally homotopic to the identity natural transformation on $\bar A$. Such a functor $\bar A$ is a \textbf{geometric realization over $\CoH_{(p)}$} of a natural coalgebra retract $A$ of $T$ 
if there is a natural isomorphism of functors from $\Top$ to the category of coalgebras 
\[
H_\ast \circ \bar{A} \cong A \circ \s^{-1}\tilde{H}_\ast. 
\]
A \textbf{natural homotopy decomposition} of $\om:\CoH_{(p)}\to \Top$ is a pair of functors $\bar A, \bar B$ from $\CoH_{(p)}\to \Top$ equipped with natural homotopy equivalences $\om \simeq \bar A \times \bar B$. 
A natural homotopy decomposition $\om\simeq \bar A \times \bar B$ is a \textbf{geometric realization over $\CoH_{(p)}$} of the natural coalgebra decomposition $T\cong A \otimes B$ if $\bar A$ and $\bar B$ are geometric realizations of $A$ and $B$, respectively.

We are now equipped to describe the functorial analogs of the decomposition methods of Cohen, Moore, and Neisendorfer. 
These results give a wonderful algebraic source of homotopy decompositions of looped co-H-spaces.

\begin{thm}{\rm\cite{STW2}}\qua 
Every natural coalgebra retract of $T$ has a geometric realization over $\CoH_{(p)}$.
\end{thm}

\begin{cor}{\rm\label{corRealizeCoalgDecomp} \cite{STW2}}\qua 
Every natural coalgebra decomposition of $T$ has a geometric realization over $\CoH_{(p)}$.
\end{cor}

%
%

We will be interested in a particular natural coalgebra decomposition of the tensor algebra functor known as the minimal decomposition, which we now set about describing. Beginning with F. Cohen, there was an interest in studying the minimal functorial coalgebra retract $\amin$ of $T$ for which $V\subseteq \amin(V)$ for every vector space $V$; 
we should note that constructions of $\amin$ are theoretical, and concrete information regarding this functor can be difficult to come by \cite{SWAmin}. 
Cohen conjectured that the primitives of $T(V)$, considered as a Hopf algebra, having tensor length not a power of $p$ must lie in the coalgebra complement of $\amin(V)$ in $T(V)$. This was confirmed by Selick and Wu, who discovered the minimal decomposition and began studying its structural properties in \cite{SelickWu}.

\begin{thm}{\rm\cite{SelickWu}\label{thmMinimalDecomp}}\qua
There is a natural coalgebra-split sub-Hopf algebra $\bmax$ of $T$ and a natural coalgebra decomposition
\begin{equation}\label{eqMinimalDecomp}
T\cong \amin\otimes\bmax. 
\end{equation} 
Moreover, $L_n(V)\subseteq \bmax(V)$ if $n$ is not a power of $p$. Here $L_n(V)$ denotes the submodule of homogeneous Lie elements of tensor length $n$ in $T(V)$. The natural coalgebra decomposition \eqref{eqMinimalDecomp} is known as the \textbf{minimal decomposition}. 
\end{thm}

By Corollary \ref{corRealizeCoalgDecomp} the minimal decomposition has a geometric realization \newline{} $\om \simeq {\bamin}\times {\bbmax}$ over $\CoH_{(p)}$. We can find more structure in this homotopy decomposition of $\om$ by making use of the observation that $\bmax$ is a natural sub-Hopf algebra of $T$. 
For a Hopf algebra $M$, write $IM$ for the augmentation ideal of $M$, and write $QM=IM/(IM)^2$ for the module of indecomposables of $M$. 
Suppose $B$ is any natural coalgebra-split sub-Hopf algebra $B$ of $T$. Since for each vector space $V$, $B(V)$ is a sub-Hopf algebra of $T(V)$, it follows that $B(V)$ is also a tensor algebra. That is, there is a natural isomorphism of algebras 
\[
B(V)\cong T\left(\bigoplus_{n\geq 1} Q_nB(V)\right) 
\]
where $Q_nB(V)$ is the image of submodule 
\[
B_n(V)=IB(V)\cap V^{\otimes n}\subseteq T(V)
\]
of $B(V)$ consisting of elements of tensor length $n$ in $T(V)$ lying in the augmentation ideal of $B(V)$ under the natural map $B(V) \to QB(V)$. The construction of each $Q_nB(V)$ is natural, and so we obtain natural isomorphisms 
\[
B \cong T \circ \bigoplus_{n\geq 1} Q_nB. 
\]
Ideally one can geometrically realize this additional structure as well; this is the content of the following theorem of Grbi\`{c}, Theriault, and Wu.

\begin{thm}{\rm\cite{DecompOmCoH}\label{thmRealizeSubHopf}}\qua
Suppose $B$ is a natural coalgebra-split sub-Hopf algebra of $T$. There exist functors $\bar{Q}_nB:\CoH_{(p)}\to \Top$ with 
\begin{enumerate}[(1)]
\item $\s^{-1} \tilde{H}_\ast(\bar{Q}_nB(Y))\cong Q_nB(\s^{-1}\tilde{H}_\ast(Y))$, 

\item $\bar{Q}_nB(Y)$ is naturally a retract of an $(n-1)$-fold desuspension of $Y^{\wedge n}$, the $n$th smash power of $Y$,

\item $\bar{B}(Y)\simeq \om\left(\bigvee_{n\geq 1} \bar{Q}_nB(Y)\right)$. 

\end{enumerate}
\end{thm}

The statement (2) requires some justification. 
 Theriault \cite{TheriaultOnProducts} has shown that if $X$ and $Y$ are coassociative co-H-spaces then $X\wedge Y\simeq \s Z$ for some co-H-space $Z$. In \cite{grayOnGenWHPs}, Gray showed that the coassociativity requirement could be relaxed -- we need only require that one of the factors in the smash product be simply-connected or a suspension space. 
 Inductively, it follows that an $n$-fold smash product of simply-connected co-H-spaces is an $(n-1)$-fold suspension of a co-H-space; symbolically, for simply-connected co-H-spaces $X_i,i=1,\dots,n$ 
\begin{equation}\label{desuspIdent}
\bigwedge_{i=1}^n X_i \simeq \s^{n-1} Z
\end{equation}
for some co-H-space $Z$. Of course there may be many choices for the space $Z$. For example, the well known decomposition 
\[
\s (X \times Y)\simeq \s X \vee \s Y \vee \s (X\wedge Y)
\]
and the failure of the identity 
\[
X\times Y \simeq X \vee Y \vee (X\wedge Y)
\]
witnesses the failure of a cancellation property for $\s$. This ambiguity need not worry us, since we will only have a need to describe the homology of a space $Z$ fitting in $\s^{n-1} Z \simeq Y^{\wedge n}$. That the space $Z$ can be chosen to admit a co-H-structure also illustrates that $\bar{Q}_nB(Y)$ can be endowed with the structure of a co-H-space, which will be of importance in the proof of Theorem \ref{thmMain}. 


%
%


\section{Decomposition Methods in Phantom Map Theory} \label{secDecompMethods}

%

In this section we develop tools which will be used to bridge the gap between the decompositions of Section \ref{secDecompCoH} and phantom map theory. 
The Loop- and Wedge-Splitting theorems (and their duals) have many applications outside our present scope, due to the existence of a vast library of decompositions in the literature to which these theorems can be applied. To substantiate this claim, we provide an application of the Loop-Splitting theorem to special cases of Question \ref{quePhantomsToFiniteCat} in Example \ref{exPhantToGanea}.
%
%
%
%
%


\begin{thm}[Loop-Splitting Theorem]
Suppose $Y$ has finite type over $\ZZ$ or $\ZZ_{(p)}$ for some prime $p$, and $\om Y\simeq A \times \om B$. If $B$ is the target of essential phantom maps from finite type domains, then so is $Y$. 
\end{thm}

%
%
%
%

\begin{proof}
Take $X$ to be an arbitrary finite type domain and write 
\[
G_n= [X,\om Y^{(n)}] \hsand H_n = [X,\om B^{(n)}].
\]
We make use of the identification
\[
\Ph(X,Y)\cong\limone G_n \hsand \Ph(X,B)\cong\limone H_n. 
\]
By Theorem \ref{limoneCtble} if $\Ph(X,Y)=\ast$ then $\{G_n\}$ is Mittag-Leffler. Since $\om Y \simeq A \times \om B$ we have a natural projection $f:\om Y \to \om B$ inducing surjections $f_n:G_n\to H_n$ of pointed sets. 

If we knew each $f_n$ was a homomorphism of groups, we could conclude $\Ph(X,B)\cong \limone H_n=\ast$ by noting $\limone f:\limone G_n\to\limone H_n$ is surjective and $\limone G_n=\ast$. In general, however, we cannot expect the functions $f_n$ to be homomorphisms, and so we must work marginally harder. 

Fortunately, the Mittag-Leffler condition makes no reference to the group structure of the individual stages of a tower, and is more a property of the underlying tower of sets. In light of Theorem \ref{limoneCtble}, to show $\limone H_n=\ast$ it suffices to show the Mittag-Leffler condition is preserved under epimorphisms of towers of pointed sets. This is the content of the following lemma. 

\begin{lem}
If $f:\{G_n\}\to\{H_n\}$ is an epimorphism of towers of pointed sets, and $\{G_n\}$ satisfies the Mittag-Leffler condition, then so does $\{H_n\}$.
\end{lem}

\noindent\textit{Proof of Lemma.} Since $\{G_n\}$ is Mittag-Leffler then for each $k$ there is some $N\in\NN$ so that for $n\geq N$ one has 
\[
G_k^{(N)}=G_k^{(n)}. 
\]
A quick diagram chase shows that the surjections $f_k:G_k\to H_k$ induce surjections $f_k^{(n)}: G_k^{(n)}\to H_k^{(n)}$. 
In other words, 
\[
H_k^{(n)} = \{ f(x) \mid x \in G_k^{(n)}\}. 
\]
But, for $n\geq N$ we have $G_k^{(n)}=G_k^{(N)}$ and so this shows $H_k^{(n)}=H_k^{(N)}$. So, the tower $\{H_n\}$ is Mittag-Leffler, which completes the proof of the lemma, and hence the proof of the Loop-Splitting Theorem. 
\end{proof}

\begin{thm}[Wedge-Splitting Theorem]
Suppose $Y$ is simply-connected and has finite type over $\ZZ$ or $\ZZ_{(p)}$ and $Y\simeq A \vee B$. If both $A$ and $B$ are rationally nontrivial, then $Y$ is the target of essential phantom maps from finite type domains. 
\end{thm}

For the proof we will need the following variation of Iriye's Corollary 1.5 from \cite{iriyePhantomsToSuspensions}. The proof is a simple modification of Iriye's argument in \cite{iriyePhantomsToSuspensions}, replacing Iriye's Theorem 2.1 with our Proposition \ref{propPLocalNoPhantoms}. 

\begin{lem}\label{lemCor1.5}
Suppose $Y$ has finite type over $\ZZ_{(p)}$. If either
\begin{enumerate}[(1)]
\item there is some $\alpha\in\pi_{2n+1}(Y)$ of infinite order whose image under the Hurewicz map is also of infinite order, or 

\item there is some $v\in H^{2n}(Y;\ZZ)$ of infinite order whose square $v^2$ is also of infinite order,
\end{enumerate}
then $\s Y$ is the target of essential phantom maps from finite type domains. 
\end{lem}

\begin{proof}[Proof of the Wedge-Splitting Theorem]
We note that since $Y$ is simply-connected, so too are $A$ and $B$. In the long fiber sequence induced by the inclusion $i:A\vee B \to A\times B$
\[
\dots \too \om F \lra{\om f} \om(A\vee B) \lra{\om i} \om A \times \om  B \lra \partial F \lra f A \vee B \lra i A \times B
\]
we can identify $F\simeq (\om A) \ast (\om B)$, where $X\ast Y$ denotes the join of topological spaces $X$ and $Y$, and we find that $\partial\simeq \ast$. It follows that $\om i$ has a section, and $\om f$ has a retraction, which gives a natural homotopy equivalence
\begin{equation}\label{eqLoopWedgeSplit}
\om(A\vee B)\simeq \om A \times \om B \times \om((\om A)\ast(\om B)). 
\end{equation}
For a more complete account of this discussion we refer the reader to the work of Porter \cite{porter}. 
We now proceed by cases. 

\vspace{.5pc}\noindent\textbf{Case I }  Suppose $Y$ has finite type over $\ZZ$. Then so do $A$ and $B$. Now, if both $A$ and $B$ are rationally nontrivial, then $(\om A)\ast(\om B)$ is a simply-connected, rationally nontrivial suspension space, hence is the target of essential phantom maps from finite type domains by Theorem \ref{thmIriye}. Applying the Loop-Splitting Theorem to the splitting \eqref{eqLoopWedgeSplit} then implies $A\vee B$ is the target of essential phantom maps from finite type domains. 

\vspace{.5pc}\noindent\textbf{Case II } In case $Y$ has finite type over $\ZZ_{(p)}$ our goal will be, as above, to show that $\om A \ast \om B$ is the target of essential phantom maps from finite type domains and appeal to the Loop-Splitting Theorem. But, since $\om A \ast \om B$ is not of finite type over $\ZZ$ we must make use of Lemma \ref{lemCor1.5}. To do so we need to discover more about $\om A \wedge \om B$. Suppose $\conn_\QQ(A)=n$ and $\conn_\QQ(B)=m$, where by $\conn_\QQ(X)=k-1$ we mean $\pi_i(X)\otimes\QQ=0$ for $i<k$ and $\pi_k(X)\otimes\QQ\neq0$. Choose $a\in H^n(\om A;\ZZ), b\in H^m(\om B;\ZZ)$ of infinite order. We proceed by cases. 

\vspace{.5pc}\noindent\textbf{Case A } If $n$ and $m$ are both even, then $a^2,b^2$ can be seen to be of infinite order, since $H^\ast(\om A;\QQ)$ contains $\QQ[\bar a]$ as a subalgebra, where $\bar a $ is the image of $a$ under rationalization, and similarly $\QQ[\bar b]$ is a subalgebra of $H^\ast(\om B;\QQ)$. Then $(a\otimes b)^2$ has infinite order in $H^\ast(\om A\wedge \om B;\ZZ)$, since $(\bar a\otimes \bar b)^2$ is nonzero in $H^\ast(\om A \wedge \om B;\QQ)$ and part (2) of Lemma \ref{lemCor1.5} applies. Here we use the K\"{u}nneth theorem to embed $H^\ast(\om A;\ZZ)\otimes H^\ast(\om B;\ZZ)$ in $H^\ast(\om A \wedge \om B;\ZZ)$ as a submodule. 

\vspace{.5pc}\noindent\textbf{Case B } If $n$ is even and $m$ is odd, then $\conn_\QQ(\om A \wedge \om B)=n+m-1$ and by the Hurewicz theorem 
$\pi_{n+m}(\om A \wedge \om B)\to H_{n+m}(\om A \wedge \om B)$ is an isomorphism, with $n+m$ odd, so part (1) of Lemma \ref{lemCor1.5} applies. 

\vspace{.5pc}\noindent\textbf{Case C } Suppose $n$ and $m$ are both odd, and without loss of generality assume $n\leq m$. Since $\conn_\QQ(\om A \wedge \om B)=n+m-1$ the rational Hurewicz homomorphism $\pi_{2n+m}\otimes \QQ\to H_{2n+m}(-;\QQ)$ is an isomorphism by the rational Hurewicz theorem. Since $n$ and $m$ are odd, $2n+m$ is odd, while $\pi_{2n+m}(\om A \wedge \om B)\otimes \QQ\neq 0$, and so part (1) of Lemma \ref{lemCor1.5} applies. 
\end{proof}

\section{Proof of Theorem \ref{thmMain}} \label{secCoH}

We begin by showing it suffices to prove Theorem \ref{thmMain} in case the nilpotent co-H-space $Y$ in question is simply-connected, so that we may appeal to the decompositions of looped co-H-spaces described in Section \ref{secDecompCoH}. 
To this end, assume $Y$ is a co-H-space with $H_i(Y;\QQ)\neq 0$ for some $i>1$. By Fox \cite{fox} $\tilde Y$ is a co-H-space, and as a consequence of the work of Iwase, Saito, and Toshio \cite{IwaseSaitoToshio} on homology of universal covers of co-H-spaces we see that if $H_i(Y;\QQ)\neq 0$ then $H_i(\tilde Y;\QQ)\neq 0$. In light of these facts and the upcoming Lemma \ref{lemSConn} we replace $Y$ with its universal cover for the proof of Theorem \ref{thmMain}. 

\begin{lem}\label{lemSConn}
Suppose $Y$ is a nilpotent co-H-space and let $c:\tilde Y \to Y$ be the universal cover. If $\tilde Y$ is the target of essential phantom maps from finite type domains, then so too is $Y$. 
\end{lem}

\begin{proof}
By Theorem \ref{thmNoPhantoms} if $\tilde Y$ is the target of essential phantom maps from finite type domains, then $\Ph(K(\ZZ,n),\tilde{Y})\neq\ast$ for some $n\geq 2$. 
We argue that $c$ induces a weak injection $\Ph(K(\ZZ,n),\tilde Y)\to \Ph(K(\ZZ,n),Y)$. 


Suppose $\phi:K(\ZZ,n)\to \tilde Y$ is an essential phantom map.
The map $c$ is the fiber of the classifying map $Y\to B\pi_1(Y)$. Since $Y$ is a co-H-space $\pi_1(Y)$ is a free group, and since $Y$ is nilpotent $\pi_1(Y)$ is either trivial or congruent to $\ZZ$. Since the result is trivial in case $\pi_1(Y)=1$ we assume $\pi_1(Y)\cong \ZZ$. So $B\pi_1(Y)\simeq S^1$ and we have a fiber sequence 
\[
\om S^1\lra\delta \tilde Y \lra c Y.
\]
We proceed by contradiction. Suppose $c\phi\simeq \ast$. Then there is a lift $\lambda:K(\ZZ,n)\to \om S^1$ of $\phi$ through $\delta$. But $\om S^1\simeq \ZZ$ is discrete and $K(\ZZ,n)$ is connected so $\lambda\simeq \ast$ and $\phi\simeq\delta\lambda$ is trivial, a contradiction. Hence $c\phi:K(\ZZ,n)\to Y$ is essential.
%
\end{proof}

We now derive Theorem \ref{thmMain} as a consequence of the following three propositions. We begin with the case $\dim_\QQ\tilde H_\ast(Y;\QQ)\geq 2$. This condition ensures the decompositions of Section \ref{secDecompCoH} are algebraically rich enough to detect essential phantom maps into $Y$ via techniques developed in Section \ref{secDecompMethods}.

%
%
%
%

\begin{prop}\label{propDimAtLeast2}
Suppose $Y$ is a simply-connected co-H-space with  $\dim_\QQ \tilde{H}_\ast(Y;\QQ)\geq 2$. Then $Y$ is the target of essential phantom maps from finite type domains. 
\end{prop}

\begin{proof}
Choose a homogeneous basis of integral classes $\{x_1,x_2,\dots\}$ for $\tilde{H}_\ast(Y;\QQ)$ with $|x_i|\leq |x_{i+1}|$ for each $i$, where $|x|$ denotes the homogeneous degree of $x$ in $\tilde{H}_\ast(Y;\QQ)$. Write  
\[
a=\s^{-1} x_1 \in \s^{-1}\tilde{H}_{m+1}(Y;\QQ) \hsand b=\s^{-1} x_2 \in \s^{-1}\tilde{H}_{n+1}(Y;\QQ).
\]
Choose a prime $p\geq 5$ so that 
\[
H_{\leq m+n+2}\left(Y^{\wedge 2};\ZZ\right) \hsand H_{\leq 2m+n+3}\left(Y^{\wedge 3};\ZZ\right)
\]
have no $p$-torsion. 
We identify $a$ and $b$ as elements of ${H}_m(\om Y;\QQ)$ and ${H}_{n}(\om Y;\QQ)$, respectively, via that Bott-Samelson Theorem. 
We will also write $a,b\in {H}_\ast(\om Y;\ZZ)$ for lifts of $a$ and $b$, and we will use the same notation for the mod $p$ reductions of these elements in ${H}_\ast(\om Y;\FF_p)$, making the context clear by indicating coefficient rings. 
We replace $Y$ with its $p$-localization to avoid cumbersome notation; that is, we write $Y$ for $Y_{(p)}$. 

To show $Y$ is the target of essential phantom maps from finite type domains, we consider the geometric realization 
\[
\om Y\simeq \om \bamin(Y)\times \om\left(\bigvee_{n\geq 2} \bar{Q}_n\bmax(Y)\right)
\]
of the minimal decomposition from Section \ref{secDecompCoH}. We justify the indexing $n\geq 2$ by noting that $Q_1\bmax=0$ since $V\subseteq \amin(V)$ for all vector spaces $V$. 
By the Loop-Splitting Theorem, it suffices to show that $\bigvee_{n\geq 2} \bar{Q}_n\bmax(Y)$ is the target of essential phantom maps from finite type domains. By the Wedge-Splitting Theorem, this will follow if $\bar{Q}_i\bmax(Y)$ is rationally nontrivial for at least two $i$. 
We set about showing this is the case. 

Write $V=\s^{-1}\tilde{H}_\ast(Y;\FF_p)$ and identify 
\[
H_\ast(\om Y;\FF_p)\cong T(V)
\]
through the Bott-Samelson Theorem. According to Theorem \ref{thmMinimalDecomp}, when $i$ is not a power of $p$ one has $L_i(V)\subseteq \bmax(V)$. So, since $p\geq 5$ we see that $[a,b],[[b,a],a]\in\bmax(V)$. Moreover, $[a,b]$ is indecomposable in $\bmax(V)$, since the tensor-length of $[a,b]$ in $T(V)$ is two, and $\bmax(V)$ contains no elements of tensor-length one in $T(V)$ (again, since $V\subseteq \amin(V)$). Similarly, $[[b,a],a]$ is indecomposable, and we have 
$[a,b] \in Q_2\bmax(V)$ and $[[b,a],a]\in Q_3\bmax(V)$. 

It follows that $[a,b]$ is in the image of 
\[
H_{n+m}\left(\om \bar{Q}_2\bmax(Y);\FF_p\right) \too H_{n+m}\left(\om Y;\FF_p\right),
\]
and so $H_{n+m}\left(\om \bar{Q}_2\bmax(Y);\FF_p\right)\neq 0$. Finally, we note $\bar{Q}_2\bmax(Y)$ is a co-H-space by Theorem \ref{thmRealizeSubHopf} and so by the Bott-Samelson Theorem
\[
H_\ast\left(\om \bar{Q}_2\bmax(Y);\FF_p\right)\cong T\left(\s^{-1}\tilde{H}_\ast\left(\bar{Q}_2\bmax(Y);\FF_p\right)\right). 
\]
Hence we infer 
\begin{equation}\label{eqQ2Fp}
\tilde{H}_{\leq m+n+1}\left(\bar{Q}_2\bmax(Y);\FF_p\right)\neq 0.
\end{equation}
 Similarly, 
 \begin{equation}\label{eqQ3Fp}
 \tilde{H}_{\leq 2m+n+1}\left(\bar{Q}_3\bmax(Y);\FF_p\right)\neq 0.
 \end{equation} 
 
Now, according to Theorem \ref{thmRealizeSubHopf} for each $i$ the space $\bar{Q}_i\bmax(Y)$ is a retract of an $(i-1)$-fold desuspension of $Y^{\wedge i}$. In particular, 
$H_k\left(\bar{Q}_i\bmax(Y);\ZZ\right)$ is a retract of  $H_{k+i-1}\left(Y^{\wedge i};\ZZ\right).$ 
So, if $H_{\leq m+n+1}\left(\bar{Q}_2\bmax(Y);\ZZ\right)$ has $p$-torsion, then so does $H_{\leq m+n+2}\left(Y^{\wedge 2};\ZZ\right)$. Similarly, if $H_{\leq 2m+n+1}\left(\bar{Q}_3\bmax(Y);\ZZ\right)$ has $p$-torsion, so does $H_{\leq 2m+n+3}\left(Y^{\wedge 3};\ZZ\right)$. 
So, since 
\[
H_{\leq m+n+2}\left(Y^{\wedge 2};\ZZ\right) \hsand H_{\leq 2m+n+3}\left(Y^{\wedge 3};\ZZ\right)
\]
have no $p$-torsion we find 
\[
\tilde{H}_{\leq m+n+2}\left(\bar{Q}_2\bmax(Y);\QQ\right) \hsand \tilde{H}_{\leq 2m+n+3}\left(\bar{Q}_3\bmax(Y);\QQ\right)
\]
are nonzero.
\end{proof}

In case $Y$ is a simply-connected finite type co-H-space with $\dim_\QQ \tilde{H}_\ast(Y;\QQ)=1$ we are unable to use the method of the proof of Proposition \ref{propDimAtLeast2} to witness the existence of essential phantom maps into $Y$ from finite type domains; we cannot expect to produce rationally nontrivial commutators in $H_\ast(\om Y;\ZZ)$, which ultimately were the driving force behind that argument. In this case $Y$ is rationally equivalent to a sphere. We proceed by cases on the parity of the dimension of this sphere. 

\begin{prop}\label{propEvenCoho}
Suppose $Y$ is a nilpotent co-H-space with $H^{2n}(Y;\QQ)\neq0$ for some $n\geq 1$. Then $Y$ is the target of essential phantom maps from finite type domains. 
\end{prop}

\begin{cor}
Suppose $Y$ is a nilpotent co-H-space with $Y\sim_\QQ S^{2n}$ for some $n\geq 1$. Then $Y$ is the target of essential phantom maps from finite type domains. 
\end{cor}

\begin{proof}[Proof of Proposition \ref{propEvenCoho}]
Let $Y\lra g K(\ZZ,2n)$ represent an element of $H^{2n}(Y;\ZZ)$ of infinite order. According to Ganea \cite{ganea1967}, since $Y$ is a co-H-space there is a lift $\lambda$ in the diagram 
\[
\xymatrix{
&& \s K(\ZZ,2n-1)\ar[d]^{p} \\
Y\ar@/^/[urr]^-{\lambda} \ar[rr]_{g} && K(\ZZ,2n), 
}
\]
where $p:\s K(\ZZ,2n-1)\simeq \s \om K(\ZZ,2n)\to K(\ZZ,2n)$ is the evaluation map.
Since $g$ induces a surjection on $\pi_{2n}\otimes\QQ$ and $p$ induces an isomorphism on $\pi_{2n}$ we can be sure $\pi_{2n}(\lambda)\otimes \QQ$ is surjective. Since $\s K(\ZZ,2n-1)$ is rationally equivalent to $S^{2n}$ we have an isomorphism of vector spaces
\begin{equation}\label{eqRatHtpySuspEML}
\pi_\ast(\s K(\ZZ,2n-1))\otimes\QQ\cong \QQ\cdot \alpha\oplus \QQ\cdot[\alpha,\alpha],
\end{equation}
where $\alpha\in \pi_{2n}(\s K(\ZZ,2n-1))\otimes\QQ$ is a nonzero element and $[-,-]$ denotes the Whitehead product. 
Since $\alpha$ is in the image of $\pi_{2n}(\lambda)$, it follows from the naturality of the Whitehead product that $\pi_\ast(\lambda)\otimes\QQ$ is surjective. 

Finally, note that by Theorem \ref{thmRatSurj} the map $\lambda: Y\to \s K(\ZZ,2n-1)$ induces surjections 
\[
\Ph(X,Y)\to \Ph(X,\s K(\ZZ,2n-1))
\]
 for all finite type spaces $X$. By Theorem \ref{thmIriye} there is a finite type space $X$ for which $\Ph(X,\s K(\ZZ,2n-1))\neq \ast$, and so $\Ph(X,Y)\neq \ast$. 
 \end{proof}

\begin{prop}\label{prop43}
If $Y$ is a nilpotent co-H-space with $Y\sim_\QQ S^{2n+1}, n\geq 1$ then $Y$ is the target of essential phantom maps from finite type domains.
\end{prop}

\begin{proof} 
We first reduce to the case where $Y$ is $(2n)$-connected.
According to Golasi\'{n}ski and Klein \cite{GolasinskiKlein} if $Y$ is a co-H-space, then one can choose compatible co-H-structures $Y$ and on each skeleton $Y_k$ so that the inclusion maps $Y_k\hookrightarrow Y$ are co-H-maps. Berstein and Hilton have shown the cofiber of a co-H-map is a co-H-space \cite[Theorem 3.4]{BersteinHilton}, so $Y/Y_k$ is a co-H-space. Finally, $Y\to Y/Y_{2n}$ is a rational equivalence, and so by Theorem \ref{thmRatSurj} this map induces a surjection $\Ph(X,Y)\to \Ph(X,Y/Y_{2n})$ for all finite type domains $X$. Hence if $Y/Y_{2n}$ is the target of essential phantom maps from finite type domains, then so too is $Y$.

Henceforth we assume the space $Y$ to be $(2n)$-connected. 
We proceed by contradiction. 
Suppose $Y$ is not the target of essential phantom maps from finite type domains. For brevity, write $K=K(\ZZ,n)$. Then by Theorem \ref{thmNoPhantoms} there is a rational equivalence $f:K\to \om Y$. 
 Let $u:\om Y\to K$ represent a cohomology class of infinite order, and write $F$ for the homotopy fiber of $u$. Since $f$ and $u$ are rational equivalences we can localize at a large enough prime $p$ and find that $f_{(p)},u_{(p)}$ induce isomorphisms on $\pi_{2n}$. For the rest of this section all spaces and maps will be localized at this large prime $p$, though the notation will not be burdened with this assumption; we write $Y$ for $Y_{(p)}$.

Now $uf$ is a self-equivalence of $K$ by the Whitehead Theorem, and so $K$ is a retract of $\om Y$. Thus $\om Y \simeq K\times F$, which gives rise to a homotopy equivalence
\[
\s \om Y \simeq \s K\vee \s F \vee \s K\wedge F. 
\]

Choose a section $s:Y\to \s \om Y$ of the evaluation map, ensured to exist since $Y$ is a co-H-space. Let $i:Y\to K$ be the composite 
\[
Y\lra s \s \om Y \simeq  \s K\vee \s F \vee \s K\wedge F \too \s K
\]
and let $q$ be the map 
\[
\s K\hookrightarrow  \s K\vee \s F \vee \s K\wedge F \simeq \s \om Y \too Y,
\]
where the last map is the evaluation map. 
Then $qi$ induces in isomorphism on $\pi_{2n+1}(Y)$. Since $Y$ is $(2n)$-connected and of finite type, it follows from the Hurewicz Theorem that
\[
q^\ast:H^{2n+1}(Y;\ZZ)\to H^{2n+1}(\s K;\ZZ)
\]
is an isomorphism. 

Now, we take a generator $v \in H^{2n+1}(\s K;\ZZ/p)$ and let $w=(q^\ast)^{-1}(v)\in H^{2n+1}(Y;\ZZ/p),$. Then $v=\s\tilde v$ for $\tilde{v}$ a generator of $H^{2n}(K;\ZZ/p)$, where $\s:H^{2n}(K)\to H^{2n+1}(\s K;\ZZ/p)$ is the suspension isomorphism. We then consider the morphism of Bockstein spectral sequences $q^\ast:E^\ast(Y)\to E^\ast(\s K)$.
Write $\P^n$ for the $n$-th reduced $p$-th power map. Since $\tilde v^p=\P^n(\tilde v)$ survives to $E^\infty_{2np}(K)$, $\P^n(v)$ survives to $E^\infty_{2np+1}(\s K)$. Since $\P^n(v)=\P^n(q^\ast(w))=q^\ast\P^n(w)$  we infer $\P^n(w)$ survives to $E^\infty_{2np+1}(Y)$. It follows that $H_{2np+1}(Y;\QQ)\neq 0$, contradicting the hypothesis $Y\sim_\QQ S^{2n+1}$. 
\end{proof}


\section{Examples} \label{secExamples}

In Examples \ref{exBHonSteroids} and \ref{ex2} we describe a co-H-spaces  meeting the hypotheses of Theorem \ref{thmMain},  but not Theorem \ref{thmIriye}. More specifically, we construct non-suspension co-H-spaces whose rational homology is nontrivial. We prefer to present infinite dimensional examples, since Zabrodsky obtained much stronger results than we have herein on phantom maps into finite complexes in the 1987 paper \cite{zabrodsky}.

\begin{ex}\label{exBHonSteroids}
For each prime $p\geq 3$ write $\alpha_p:S^{2p}\to S^3$ for a representative of an element of order $p$ in $\pi_{2p}(S^3)$. The homotopy cofibers $C_{\alpha_p}$ of these maps are classical examples, due to Berstein and Hilton \cite[pg. 444]{BersteinHilton}, of co-H-spaces that do not have the homotopy type of suspension spaces. One key to establishing these examples is to prove, via Berstein-Hilton-Hopf invariant techniques, that each map $\alpha_p$ is a co-H-map. By \cite[Theorem 3.4]{BersteinHilton} the cofiber of a co-H-map is a co-H-space.

Write $\alpha:\bigvee_{p\geq 3} S^{2p}\to S^3$, where the wedge is taken over all odd primes, for the map whose restriction to each summand $S^{2p}$ is $\alpha_p$. Since each $\alpha_p$ is a co-H-map, so is $\alpha$. It follows that the homotopy cofiber $C_\alpha$ of $\alpha$ is a co-H-space. Evidently $\dim_\QQ\tilde{H}_\ast(C_\alpha;\QQ)=\infty$.

We now argue that $C_\alpha$ is not a suspension space. Assume to the contrary that $C_\alpha\simeq \s Z$. Then by the proof of \cite[Lemma 3.6]{BersteinHilton} we can choose $Z$ 1-connected, so that $Z$ has a homology decomposition, i.e. there is a diagram 
\[
\xymatrix{
M_1\ar[d]_{k_1} & M_2\ar[d]_{k_2} & \dots & M_n \ar[d] & M_{n+1} \ar[d] & \dots \\
Z_1 \ar[r]_{i_1} & Z_2\ar[r]_{i_2} & \dots \ar[r] & Z_n\ar[r]_{i_n} & Z_{n+1} \ar[r] & \dots 
}
 \]
 in which each $M_i=M(H_{i+1}(Z),i), M_i\to Z_i\to Z_{i+1}$ is a cofiber sequence, 
and $Z$ is the homotopy colimit of the tower along the bottom of this diagram. 
The space $Z_i$ is called the $i$th stage of the homology decomposition. 
  It follows that $\s Z$ has a homology decomposition in which each stage is a suspension. 

Suppose $h:C_\alpha\to \s Z$ is a homotopy equivalence. Write $(C_\alpha)_k$ for the $k$th stage of the homology decomposition for $C_\alpha$. 
According to Arkowitz \cite[Proposition 3.4]{ArkowitzHomologyDecomp} 
since $\Ext(H_n(C_\alpha;\ZZ);H_{n+1}(\s Z;\ZZ))=0$ for all $n$ and $\s Z$ is 2-connected, $h$ induces homotopy equivalences  $h_n:(C_\alpha)_n\to (\s Z)_n\simeq \s (Z_n)$. But then 
$(C_\alpha)_6\simeq C_{\alpha_3}$ must be a suspension space, a contradiction. 
\end{ex}

\begin{ex}\label{ex2}
By modifying the construction from Example \ref{exBHonSteroids} we can obtain an infinite-dimensional, non-suspension co-H-space $Y$ with $Y\sim_\QQ S^3$. Replace each map $\alpha_p:S^{2p}\to S^3$ with a map $\beta_p:M(\ZZ/p,2p)\to S^3$ representing an element of $\pi_{2p}(S^3;\ZZ/p)$ of order $p$. The argument of Berstein and Hilton \cite{BersteinHilton} shows that the cofiber $C_{\beta_p}$ of each $\beta_p$ is a co-H-space which is not a suspension space, and so the argument in Example \ref{exBHonSteroids} shows that $C_\beta$ is a co-H-space which is not a suspension space. 
\end{ex}

Finally we present an application of the Loop-Splitting theorem to spaces that are not necessarily co-H-spaces. 
For a space $Y$ write $G_m(Y)$ for the $m$th space of Ganea over $Y$ (see 
\cite{cornealuptonopreatanre};
the reader may more readily recognize this space as $G_m(Y)=B_m \om Y$ where $B_m$ is the $m$th stage of Milnor's classifying space construction). The spaces $G_m(Y)$ can be thought of as prototypes for spaces of \LS category at most $m$. We view this example as a test case for Question \ref{quePhantomsToFiniteCat}.

\begin{ex} \label{exPhantToGanea}
We show that if $H_i(G_m(Y);\QQ)\neq 0$ for some $i>1$ then $G_m(Y)$ is the target of essential phantom maps from finite type domains. 

There is a well-known homotopy decomposition
\[
\om G_m(Y)\simeq \om Y \times \om((\om Y)^{\ast m+1})
\]
where $X^{\ast k}$ denotes the $k$-fold join of $X$. Since $H_i(G_m(Y);\QQ)\neq 0$ we must have $H_j(Y;\QQ)\neq 0$ for some $j>1$ and similarly $H_\ast((\om Y)^{\ast m+1};\QQ)$ is similarly nontrivial, so by Theorem \ref{thmIriye} $(\om Y)^{\ast m+1}$ is the target of essential phantom maps from finite type domain.  The Loop-Splitting theorem then implies $G_m(Y)$ is the target of essential phantom maps. 
\end{ex} 

\bibliography{DissertationBiblio}{}
  \bibliographystyle{plain}

\end{document}